\documentclass{svproc}
\usepackage{url}
\usepackage{graphicx,multicol,footmisc}

\usepackage[style=numeric,sorting=none]{biblatex} %Imports biblatex package
\usepackage{amsmath,amssymb}

\begin{document}
\mainmatter              % start of a contribution
\title{Bifurcations in inertial focusing of a particle suspended in flow through curved rectangular ducts}
\titlerunning{Bifurcations in inertial particle focusing}  % abbreviated title (for running head)
\author{Rahil N. Valani\inst{1} \and Brendan Harding\inst{2} \and Yvonne M. Stokes\inst{1}}
\authorrunning{Valani et al.} % abbreviated author list (for running head)
\tocauthor{Rahil N. Valani, Brendan Harding and Yvonne M. Stokes}
\institute{School of Mathematical Sciences, University of Adelaide, South Australia 5005, Australia\\
\email{rahil.valani@adelaide.edu.au},
%\\ WWW home page:
%\texttt{http://users/\homedir iekeland/web/welcome.html}
\and
School of Mathematics and Statistics, Victoria University of Wellington, Wellington 6012, New Zealand}

\maketitle              % typeset the title of the contribution

\begin{abstract}
Particles suspended in a fluid flow through a curved duct can focus to specific locations within the duct cross-section. This particle focusing is a result of a balance between two dominant forces acting on the particle: (i) the inertial lift force arising from small but non-negligible inertia of the fluid, and (ii) the secondary drag force due to the cross-sectional vortices induced by the curvature of the duct. By adopting a simplified particle dynamics model developed by Ha et al.~\cite{Kyung2021}, we investigate both analytically and numerically, the particle equilibria and their bifurcations when a small particle is suspended in low-flow-rate fluid flow through a curved duct having a $2\times1$ and a $1\times2$ rectangular cross-section. In certain parameter regimes of the model, we analytically obtain the particle equilibria and deduce their stability, while for other parameter regimes, we numerically calculate the particle equilibria and stability. Moreover, we observe a number of different bifurcations in particle equilibria such as saddle-node, pitchfork and Hopf, as the model parameters are varied. These results may aid in the design of inertial microfluidic devices aimed at particle separation by size.
% We would like to encourage you to list your keywords within
% the abstract section using the \keywords{...} command.
\keywords{Bifurcations $\cdot$ Inertial lift $\cdot$ Inertial focusing $\cdot$ Inertial particle focusing $\cdot$ Inertial microfluidics $\cdot$ Particle separation}
\end{abstract}
\section{Introduction}
Dynamics of a particle suspended in a fluid flow is governed by hydrodynamic forces
acting on the particle from the surrounding flow. At relatively low Reynolds number, where the inertia of the fluid flow is small but non-negligible, inertial lift force acts on the particle and facilitates migration of the particle across streamlines of the background flow. This was first demonstrated in the classical experiment of Segré \& Silberberg~\cite{SEGRE1961} where particles suspended in flow through a straight pipe with a circular cross-section were observed to migrate to an annular region located at approximately $0.6$ times the radius of the pipe from the pipe center. The phenomenon of inertial migration has found many applications in medical and industrial settings such as isolation of circulating tumor cells (CTCs)~\cite{CTCs3,CTCs1,CTCs2}, separation of particles and cells~\cite{Cells2,Cells1,Cells4,Cells3}, flow cytometry~\cite{flowcyto}, water filtration~\cite{water1}, extraction of blood plasma~\cite{plasma} and identification of small-scale pollutants in environmental samples~\cite{pollutant}. Recent advances in inertial microfluidics are provided in several review articles~\cite{review2,review3,Stoecklein2019,review1}.

A commonly used duct geometry in inertial microfluidics experiments is a spiral or a circular duct with a rectangular cross-section. Hence, understanding of particle focusing in curved duct geometries is crucial for designing microfluidic devices that optimize particle separation. Harding et al.~\cite{harding_stokes_bertozzi_2019} developed a general asymptotic model for forces that govern the motion of a spherical particle suspended in a fluid flow through a curved duct at low flow rates. This model was then used to investigate the inertial migration of a neutrally buoyant spherical particle suspended in flow through curved ducts with square, rectangular and trapezoidal cross-sections. They identified stable and unstable particle equilibria in the cross-section of the duct and also observed that the location and nature of these particle equilibria vary with the cross-sectional geometry, bend radius of the duct and particle size. Ha et al.~\cite{Kyung2021} developed a reduced order model, the Zero Level Fit (ZeLF) model, by fitting curves to the inertial lift force field and the secondary drag force field calculated from the simulated data of Harding et al.~\cite{harding_stokes_bertozzi_2019}. Although simulation data from \cite{harding_stokes_bertozzi_2019} can be interpolated and integrated directly to simulate particle trajectories, a simpler particle dynamics model allows for analytical treatment as well as rapid prototyping. Using this simplified model, Ha et al.~\cite{Kyung2021} numerically investigated various dynamical behaviors of small particles and the bifurcations in the particle equilibria for flow through a curved duct with a square cross-section. It has been shown that rectangular ducts are better at separating particles of different sizes compared to square ducts~\cite{harding_stokes_bertozzi_2019}. With this motivation, we build on the work of Ha et al.~\cite{Kyung2021} and explore using the ZeLF model, both analytically and numerically, the bifurcations taking place in particle equilibria for small particles in a $2\times1$ and a $1\times2$ rectangular cross-section. 

The paper is organized as follows. In Sec.~\ref{sec: theory} we briefly outline the simplified ZeLF model of Ha et al.~\cite{Kyung2021} that has been adapted in this paper for a $2\times1$ and a $1\times2$ rectangular cross-section. In Secs.~\ref{sec: 2x1_bifur} and \ref{sec: 1x2_bifur} we present the various particle equilibria and the bifurcations between them as a function of the duct bend radii for a $2\times1$ and a $1\times2$ rectangular cross-section, respectively. Finally, we provide conclusions in Sec.~\ref{sec: concl}.

\section{Theoretical Model}\label{sec: theory}

\begin{figure}
\centering
\includegraphics[width=0.95\columnwidth]{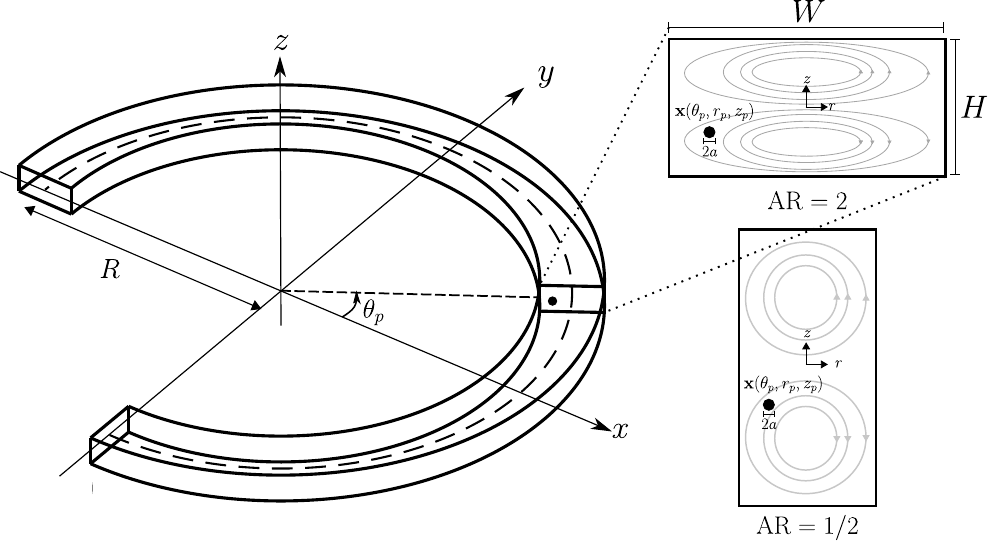}
\caption{Schematic showing the theoretical setup. A particle (black filled circle) of radius $a$ with its center located at $\mathbf{x}_p=\mathbf{x}(\theta_p,r_p,z_p)$ is suspended in a fluid flow through a curved duct of radius $R$ having a rectangular cross-section of width $W$ and height $H$ with aspect ratio defined as $\text{AR}=W/H$. Enlarged view of the two cross-sections ($\text{AR}=2$ and $\text{AR}=1/2$) show the local cross-sectional $(r,z)$ co-ordinate system, and the secondary flow (gray closed curves) induced by the curvature of the duct.}
\label{Fig: schematic}
\end{figure}

As shown in Figure~\ref{Fig: schematic}, consider a neutrally buoyant particle of radius $a$ suspended in a fluid of constant density $\rho$ and dynamic viscosity $\mu$. The fluid flows along a curved duct of constant radius $R$ having a uniform rectangular cross-sectional geometry of width $W$ and height $H$, with aspect ratio defined as $\text{AR}=W/H$. The horizontal and vertical co-ordinates within the cross-section are described using $r$ and $z$ co-ordinates with the origin located at the center of the rectangle so that the domain is $-W/2\leq r \leq W/2$ and $-H/2\leq z \leq H/2$. These cross-sectional co-ordinates are related to the global co-ordinates of the three-dimensional duct geometry as follows:
\begin{equation*}
    \mathbf{x}(\theta,r,z) = (R+r)\cos(\theta) \,{\mathbf{i}}+(R+r)\sin(\theta)\,{\mathbf{j}}+z\,{\mathbf{k}}.
\end{equation*}
Here $\theta$ is the angular co-ordinate along the curved duct and the center of the particle is located at $\mathbf{x}_p=\mathbf{x}(\theta_p,r_p,z_p)$. We now briefly outline the Zero Level Fit (ZeLF) model of Ha et al.~\cite{Kyung2021} that will be used for the analysis presented in this paper.

Based on the leading order force model of Harding et al.~\cite{harding_stokes_bertozzi_2019}, Ha et al.~\cite{Kyung2021} developed a simplified model for the dynamics of a small neutrally buoyant particle suspended in a slow flow through curved ducts with relatively large bend radii. In this model, the inertial lift force $\mathbf{L}=(L_r,L_z)$ and the secondary drag force $\mathbf{D}=(D_r,D_z)$ are approximated by fitting simple model functions to the simulated data from \cite{harding_stokes_bertozzi_2019} such that the topology of the zero level sets is preserved. The inertial lift force is approximated by that for the case of flow through a straight duct in the limit of large bend radii, and the drag force $D$ is approximated by the Stokes drag due to the secondary velocity field arising from cross-sectional vortices, again in the limit of large bend radii. Ignoring the particle inertia, the first order ODEs describing the motion of a small particle in the $r$-$z$ plane under the influence of these two force fields are given as follows:
\begin{equation}\label{eq: dim_zelf_r}
    6\pi \mu a\frac{d {r}_p}{d {t}}={{L}_r({r}_p,{z}_p)}+{D}_r({r}_p,{z}_p),
\end{equation}
\begin{equation}\label{eq: dim_zelf_z}
    6\pi \mu a\frac{d {z}_p}{d {t}}={{L}_z({r}_p,{z}_p)}+{D}_z({r}_p,{z}_p).
\end{equation}
The corresponding non-dimensional equations of motion are given by:
%This model was constructed by fitting the inertial lift $(\mathbf{L}=(L_r,L_z))$ and secondary drag $(\mathbf{D}=(D_r,D_z))$ curves from the simulated data of the Finite Element Model for a straight duct ($R\xrightarrow{}\infty$) using a combination of elementary functions in a specific way that preserves the topology of the zero level sets for the inertial lift and drag forces. In this model, the cross-sectional dynamics of the particle is governed by the following system of ODEs 
%\RV{BH-Here we use $r_p,z_p$ but in the other model we used $x_{p,r}$ and $x_{p,z}$. We should probably pick a consistent notation for this paper.}:
\begin{equation}\label{Zelf_r}
    \frac{d \tilde{r}_p}{d \tilde{t}}=\frac{1}{6\pi}\left(\frac{\tilde{a}^3}{8}{\tilde{L}_r(\tilde{r}_p,\tilde{z}_p)}+\frac{1}{2\tilde{R}}\tilde{D}_r(\tilde{r}_p,\tilde{z}_p)\right),
\end{equation}
\begin{equation}\label{Zelf_z}
    \frac{d \tilde{z}_p}{d \tilde{t}}=\frac{1}{6\pi}\left(\frac{\tilde{a}^3}{8}{\tilde{L}_z(\tilde{r}_p,\tilde{z}_p)}+\frac{1}{2\tilde{R}}\tilde{D}_z(\tilde{r}_p,\tilde{z}_p)\right).
\end{equation} 
Here the dimensionless variables denoted by tildes are defined as follows:
\begin{align*}
  (\tilde{r}_p,\tilde{z}_p)&=\frac{2}{l} (r_p,z_p),\:\:\tilde{t}=\frac{2\rho U_m^2}{\mu} t,\:\:\tilde{a}=\frac{2a}{l},\:\:\tilde{R}=\frac{2R}{l},\\ \nonumber
   (\tilde{L}_r,\tilde{L}_z)&=\frac{l^2}{\rho U_m^2 a^4}(L_r,L_z),\:\: (\tilde{D}_r,\tilde{D}_z)=\frac{4 R}{\rho a U_m^2 l^2} (D_r,D_z),\\ \nonumber 
\end{align*}
where $U_m$ is a characteristic axial flow speed of the background fluid flow (taken as the maximum axial speed) and $l=\text{min}\{W,H\}$ is a characteristic length scale of the duct cross-section. For the remainder of the paper, we will use only dimensionless variables and drop the tilde from them (except $\tilde{a}$ and $\tilde{R}$) for convenience. In this paper, we use this model to investigate the particle equilibria and the bifurcations in them for rectangular cross-sections with aspect ratios $\text{AR}=2$ and $\text{AR}=1/2$. 

\begin{figure}[t]
\centering
\includegraphics[width=\columnwidth]{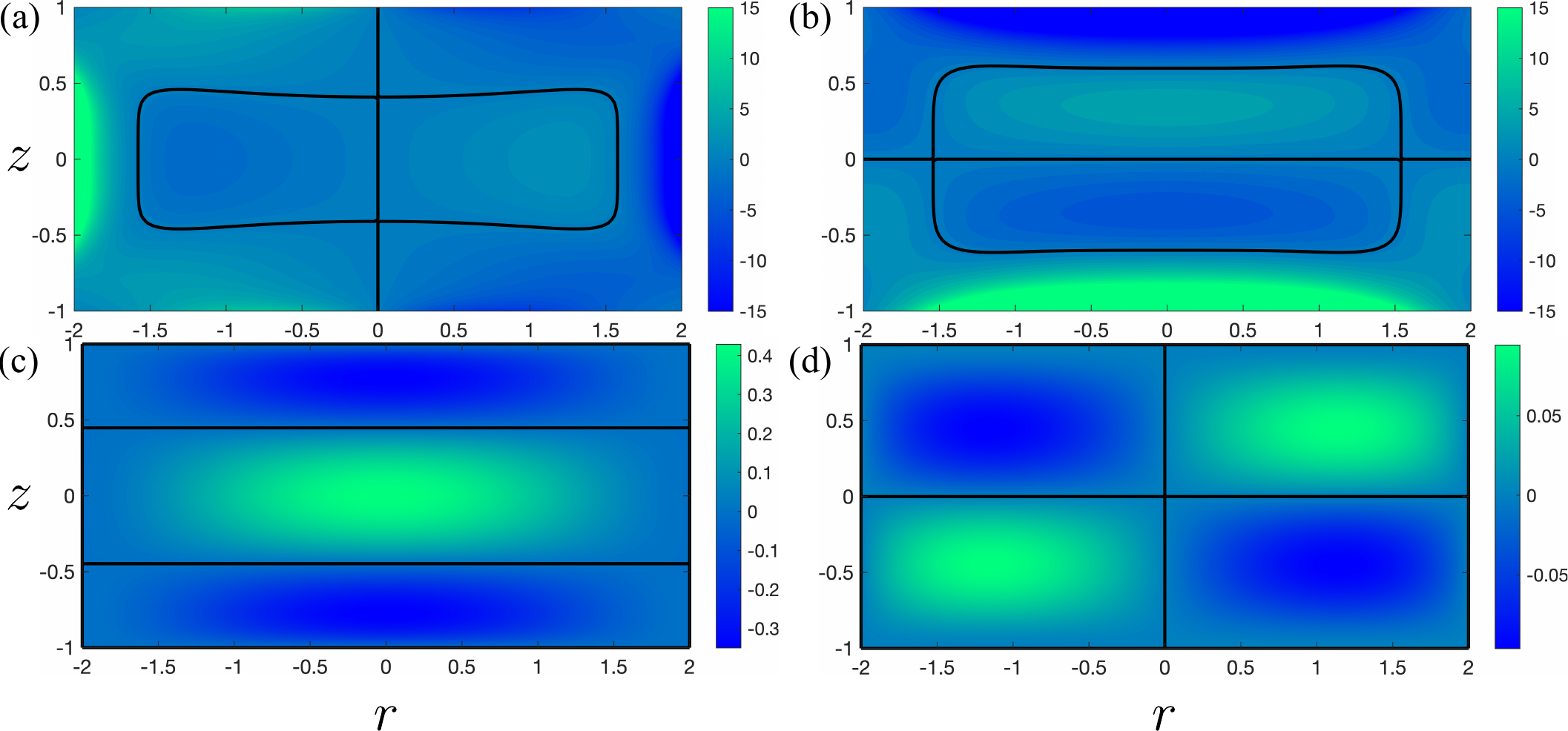}
\caption{Contour plot showing the distribution of dimensionless inertial lift and secondary drag forces inside a $2\times1$ rectangular cross-section. (a) $L^{2\times1}_r(r,z)$, (b) $L^{2\times1}_z(r,z)$, (c) $D^{2\times1}_r(r,z)$ and (d) $D^{2\times1}_z(r,z)$ are shown. The black curves in each panel show the zero level contours of the corresponding force field.}
\label{Fig: DSTA contour}
\end{figure}

\subsection{ZeLF model for a $2\times1$ rectangular cross-section}

We start by presenting the ZeLF model fits to the inertial lift force and the secondary drag force for a $2\times1$ rectangular cross-section. We fit simple model functions to the inertial lift force field from the simulated results of \cite{harding_stokes_bertozzi_2019} as was done in \cite{Kyung2021}. Doing this for a small particle of radius $\tilde{a}=0.05$ in a straight duct with a $2\times1$ rectangular cross-section one obtains:

\begin{align*}
L^{2\times1}_r (r,z) &=r\left[1-0.0643 r^6 - 25.5128 z^6 -31.1(1-0.4006 r^2)z^4\right]\\ \nonumber
& \text{exp}\Big( 0.505 + 0.427 r^2 - 5.081 z^2 - 0.2 r^4 +1.518 r^2 z^2 +0.594 z^4  \\ \nonumber
& +0.042 r^6+0.007 r^4 z^2 - 2.283 r^2 z^4 + 2.8 z^6 \Big),\\\nonumber
L^{2\times1}_z (r,z) &= z\left[1-9.0878 z^6 - 0.0316 r^8 - 1.6(1-0.1778 r^4)z^2 \right] \\ \nonumber
& \text{exp}\Big( 3.030 - 1.168 z^2 - 0.536 r^2 - 2.199 z^4 +0.476 z^2 r^2 + 0.104 r^4  \\ \nonumber
& +2.094 z^6 +0.051 z^4 r^2 - 0.212 z^2 r^4 -0.033 r^6 \Big).\\\nonumber
%D_r (r_p,z_p) &= C (1-0.25 r_p^2)^2(1-z_p^2)(1-5z_p^2)\\\nonumber
%D_z (r_p,z_p) &= C r_p z_p(1-0.25 r_p^2)(1-z_p^2)^2\\\nonumber
\end{align*}
Similarly, the fitting of the secondary drag force field yields the following:
\begin{align*}
D^{2\times1}_r (r,z) &= 6 \pi C (1-0.25 r^2)^2(1-z^2)(1-5z^2),\\\nonumber
D^{2\times1}_z (r,z) &= 6 \pi C r z(1-0.25 r^2)(1-z^2)^2.\\\nonumber
\end{align*}
where $C=0.02319$. Contours of these force fields inside a $2\times1$ rectangular cross-section are shown in Fig.~\ref{Fig: DSTA contour}. These fits compared with the simulated data of \cite{harding_stokes_bertozzi_2019} have a relative error of $9\%$ for $L^{2\times1}_r$, $3\%$ for $L^{2\times1}_z$ and $10\%$ for $D^{2\times1}_r$ and $D^{2\times1}_z$. 

\subsection{ZeLF model for a $1\times2$ rectangular cross-section}\label{Sec: zelf 1x2}

For a $1\times2$ rectangular cross-section, the inertial lift force field is simply approximated by a $90^{\circ}$ rotation of the inertial lift force field for a $2\times1$ rectangular cross-section i.e. we replace the $r$ coordinates by the $z$ coordinates and vice versa. This results in the following:  
\begin{align*}
L^{1\times2}_r (r_p,z_p) = L^{2\times1}_z (z_p,r_p), \\ \nonumber
L^{1\times2}_z (r_p,z_p) = L^{2\times1}_r (z_p,r_p). \\ \nonumber
\end{align*}
Fitting the secondary drag force field yields
\begin{align*}
D^{1\times2}_r (r,z) &= 6 \pi C (1-r^2)^2(1-0.25 z^2)(1-1.25 z^2),\\\nonumber
D^{1\times2}_z (r,z) &= 24 \pi C r z(1- r^2)(1-0.25 z^2)^2.\\\nonumber
\end{align*}
\begin{figure}[t]
\centering
\includegraphics[width=\columnwidth]{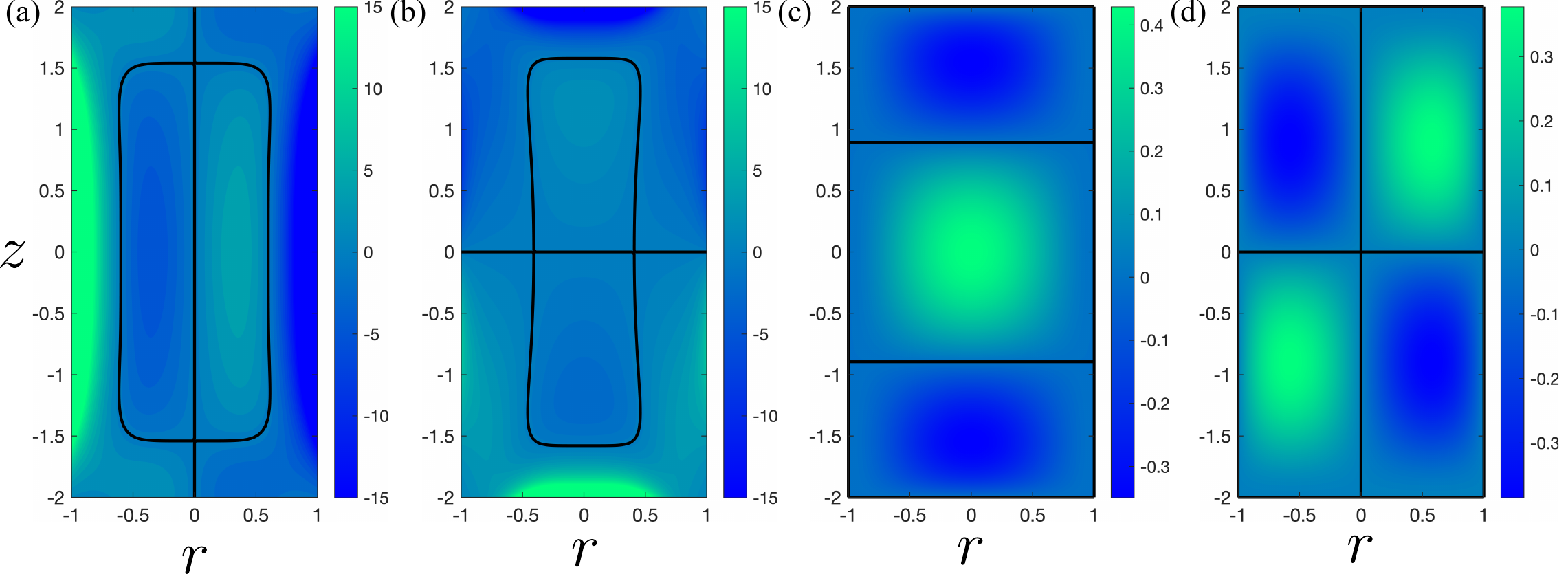}
\caption{Contour plot showing the distribution of inertial lift and secondary drag forces inside a $1\times2$ rectangular cross-section. (a) $L^{1\times2}_r(r,z)$, (b) $L^{1\times2}_z(r,z)$, (c) $D^{1\times2}_r(r,z)$ and (d) $D^{1\times2}_z(r,z)$ are shown. The black curves in each panel show the zero level contours of the corresponding force field.}
\label{Fig: DSTA contour 1x2}
\end{figure}
where $C=0.018542$. Contours of these force fields inside the $1\times2$ rectangular cross-section are depicted in Fig.~\ref{Fig: DSTA contour 1x2}. These fits for drag forces compared with the simulated data of \cite{harding_stokes_bertozzi_2019} have a relative error of $14\%$. 

\subsection{Numerical solution and the stability of particle equilibria}\label{sec: numer}

%The appropriate fitting results in simple model functions for the inertial lift $(\mathbf{L}=(L_r,L_z))$ and secondary drag $(\mathbf{D}=(D_r,D_z))$ depending on the cross-sectional shape of the duct.
%
%\section{Particle equilibria and bifurcations in a $2\times1$ rectangular duct}\label{sec: bifur}
%
The particle equilibria $(r_p,z_p)=(r^*,z^*)$, that correspond to the fixed points of the dynamical system, can be obtained from Eqs.~\eqref{Zelf_r} and \eqref{Zelf_z} by setting the time derivatives to zero giving us the following nonlinear equations to solve:
\begin{align*}
  &\frac{\tilde{a}^3}{8}{L}_r(r^*,z^*)+\frac{1}{2\tilde{R}}{D}_r(r^*,z^*)=0,\\ \nonumber
  &\frac{\tilde{a}^3}{8}{L}_z(r^*,z^*)+\frac{1}{2\tilde{R}}{D}_z(r^*,z^*)=0.
\end{align*}
To characterize the stability of particle equilibria, one can apply a small perturbation to the equilibrium state according to $(r_p,z_p)=(r^*,z^*)+\epsilon(r_1,z_1)$, where $\epsilon>0$ is a small perturbation parameter. Substituting this in Eqs.~\eqref{Zelf_r} and \eqref{Zelf_z} and comparing $\mathcal{O}(\epsilon)$ terms one gets the following linear matrix equation for the evolution of the perturbations:
\begin{gather}\label{lin_sys_2x1}
 \begin{bmatrix} \dot{r}_1 \\ \dot{z}_1 \end{bmatrix}
 = \frac{1}{6 \pi}
  \begin{bmatrix}
  \frac{\tilde{a}^3}{8}\frac{\partial L_r}{\partial r_p} + \frac{1}{2\tilde{R}}\frac{\partial D_r}{\partial r_p} &
\frac{\tilde{a}^3}{8}\frac{\partial L_r}{\partial z_p} + \frac{1}{2\tilde{R}}\frac{\partial D_r}{\partial z_p} \\
 \frac{\tilde{a}^3}{8} \frac{\partial L_z}{\partial r_p} +\frac{1}{2\tilde{R}} \frac{\partial D_z}{\partial r_p} &
  \frac{\tilde{a}^3}{8}\frac{\partial L_z}{\partial z_p} +\frac{1}{2\tilde{R}}\frac{\partial D_z}{\partial z_p}
   \end{bmatrix}_{(r^*,z^*)}
    \begin{bmatrix} {r}_1 \\ {z}_1 \end{bmatrix}.
\end{gather}
The nature of a particle equilibrium is determined by the eigenvalues $\lambda$ of the matrix on the right hand side of Eq.~\eqref{lin_sys_2x1}. The particle trajectories presented in this paper are obtained by numerically solving Eqs.~\eqref{Zelf_r} and \eqref{Zelf_z} in \texttt{MATLAB} using the inbuilt ode45 and ode15s solvers.

%and corresponding bifurcations resulting from the solving the above nonlinear set of ODEs for a rectangular $2\times1$ cross-section are shown in Figure~\ref{Fig: ZeLF model}(b). Comparing the bifurcation from the Finite Element Model in Figure~\ref{Fig: rect_2x1} with the ZeLF model in Figure~\ref{Fig: ZeLF model}(b), we again see a good agreement in the bifurcation expect for a narrow region where the stable and unstable nodes merge. In the ZeLF model we find the emergence of an unstable spiral first which changes into a stable spiral while in the FEM model, we directly observe the stable spiral. Other comparisons are similar to what was observed for the square cross-section. Similar to the square cross-section, we analytically investigate the linear stability and the fixed points in the limit of large and small $\kappa$.

\section{Particle equilibria and bifurcations in a $2\times1$ rectangular cross-section} \label{sec: 2x1_bifur}

\subsection{Large $\tilde{R}$ limit}

For large bend radius $\tilde{R}$, the inertial lift force dominates the secondary drag force on the particle. Thus, in the limit $\tilde{R}\xrightarrow{}\infty$, the secondary drag vanishes and the nonlinear system of ODEs in Eqs.~\eqref{Zelf_r} and \eqref{Zelf_z} reduces to
\begin{equation}\label{small_kappa_2x1}
    \frac{d r_p}{d t}=\frac{\tilde{a}^3}{48 \pi}{L^{2\times1}_r(r_p,z_p)}\:\:\:\text{and}\:\:\:\frac{d z_p}{d t}=\frac{\tilde{a}^3}{48 \pi}{L^{2\times1}_z(r_p,z_p)}.
\end{equation}
Solving for the fixed points of Eq.~\eqref{small_kappa_2x1} one can analytically obtain the following nine points that correspond to the particle equilibria: 
%(Note that we also get higher order fixed points corresponding to the boundaries of the cross-section but they are not of interest for the present analysis)
\begin{align*}
(r^*,z^*)=&\big{\{}(0,0), (0,\pm 0.6000) , (\pm 1.5800,0),\\ \nonumber
&(\pm1.5303,0.4094), (\pm 1.5303,-0.4094)\big{\}}.
\end{align*}
To determine the stability of these fixed points, we calculate the eigenvalues of the linear stability matrix, Eq.~\eqref{lin_sys_2x1}. These eigenvalues and the corresponding nature of the fixed points are summarized in Table~\ref{table:ZeLF rect_2x1}.

%\begin{table}
%\caption{Eigenvalues $(\lambda_1,\lambda_2)$ and eigenvectors $(V_1,V_2)$ for the fixed points $(r^*,z^*)$ of the ZeLF model with a rectangular $2\times1$ cross-section in the limit of small $\kappa$.} % title of Table
%\centering % used for centering table
%\begin{tabular}{c c c c c c} % centered columns (4 columns)
%\hline\hline %inserts double horizontal lines
%Location & Type & \lambda_1 & \lambda_2 & \mathbf{V}_1 & \mathbf{V}_2 \\ [0.5ex] % inserts table
%heading
%\hline % inserts single horizontal line
%\\
%(0,0) & unstable node & 0.0879 & 1.0980 & \begin{bmatrix} 1\\ 0 \end{bmatrix} & \begin{bmatrix} 0\\ 1 \end{bmatrix} \\[3ex]  % inserting body of the table
%
%(0,\pm 0.6000) & stable node & -0.0733 & -2.2100 & \begin{bmatrix} 1\\ 0 \end{bmatrix} & \begin{bmatrix} 0\\ 1 \end{bmatrix} \\[3ex] % inserting body of the table 
%
%(\pm 1.5800,0) & stable node & -0.0747 & -0.8464 & \begin{bmatrix} 0\\ 1 \end{bmatrix} & \begin{bmatrix} 1\\ 0 \end{bmatrix} \\[3ex] % inserting body of the table 
%
%(\pm 1.5304,0.4092) & saddle & -0.6161 & 0.2263 & \begin{bmatrix} 0.7381\\ \pm 0.6747 \end{bmatrix} & \begin{bmatrix} \mp 0.5218\\ 0.8531 \end{bmatrix} \\[3ex] % inserting body of the table
%
%(\pm 1.5304,-0.4092) & saddle & -0.6161 & 0.2263 & \begin{bmatrix} 0.7381\\ \mp 0.6747 \end{bmatrix} & \begin{bmatrix} \pm 0.5218\\ 0.8531 \end{bmatrix} \\[3ex] % inserting body of the table
%
%[1ex] % [1ex] adds vertical space
%\hline %inserts single line
%\end{tabular}
%\label{table:ZeLF rect_2x1} % is used to refer this table in the text
%\end{table}
\begin{figure}[t]
\centering
\includegraphics[width=\columnwidth]{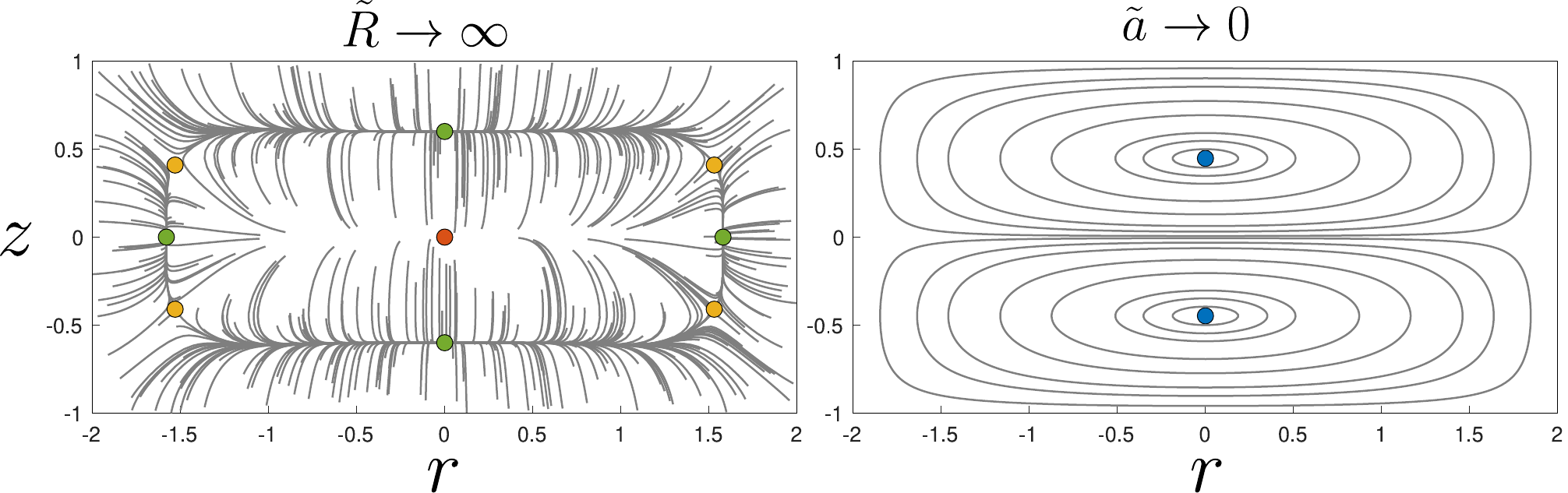}
\caption{Particle equilibria (filled circles) and particle trajectories (gray curves) in the limit of $\tilde{R} \xrightarrow{} \infty$ (left) and $\tilde{a} \xrightarrow{}0$ (right) for a $2\times1$ rectangular cross-section. The color of the filled circles indicates the type of particle equilibria obtained from linear stability analysis: unstable node in red, saddle point in yellow, stable node in green and a center in blue.}
\label{Fig: fixedpt_2x1}
\end{figure}

\begin{table}
\caption{Eigenvalues $(\lambda_1,\lambda_2)$ and eigenvectors $(\mathbf{V}_1,\mathbf{V}_2)$ for the fixed points $(r^*,z^*)$ of the ZeLF model with a rectangular $2\times1$ cross-section in the limit $\tilde{R}\xrightarrow{}\infty$.} % title of Table
\centering % used for centering table
\begin{tabular}{c c c c c c} % centered columns (4 columns)
\hline\hline %inserts double horizontal lines
Location $(r^*,z^*)$ & Type & $\lambda_1$ & $\lambda_2$ & $\mathbf{V}_1$ & $\mathbf{V}_2$ \\ [0.5ex] % inserts table
\hline % inserts single horizontal line
\\
(0,0) & unstable node & 0.0110\,$\tilde{a}^3$ & 0.1373\,$\tilde{a}^3$ & $\begin{bmatrix} 1\\ 0 \end{bmatrix}$ & $\begin{bmatrix} 0\\ 1 \end{bmatrix}$ \\[3ex]  % inserting body of the table
$(0,\pm 0.6000)$ & stable node & $-0.0092\,\tilde{a}^3$ & $-0.2762\,\tilde{a}^3$ & $\begin{bmatrix} 1\\ 0 \end{bmatrix}$ & $\begin{bmatrix} 0\\ 1 \end{bmatrix}$ \\[3ex] % inserting body of the table 
$(\pm 1.5800,0)$ & stable node & $-0.1059\,\tilde{a}^3$ & $-0.0094\,\tilde{a}^3$ & $\begin{bmatrix} 1\\ 0 \end{bmatrix}$ & $\begin{bmatrix} 0\\ 1 \end{bmatrix}$ \\[3ex] % inserting body of the table 
$(\pm 1.5304,0.4092)$ & saddle & $-0.0770\,\tilde{a}^3$ & $0.0283\,\tilde{a}^3$ & $\begin{bmatrix} 0.7380\\ \pm 0.6748 \end{bmatrix}$ & $\begin{bmatrix} \mp 0.5225\\ 0.8526 \end{bmatrix}$ \\[3ex] % inserting body of the table
($\pm$ 1.5304,-0.4092) & saddle & -0.0770\,$\tilde{a}^3$ & 0.0283\,$\tilde{a}^3$ & $\begin{bmatrix} 0.7380\\ \mp 0.6748 \end{bmatrix}$ & $\begin{bmatrix} \pm 0.5225\\ 0.8526 \end{bmatrix}$ \\[3ex] % inserting body of the table
%[1ex] % [1ex] adds vertical space
\hline %inserts single line
\end{tabular}
\label{table:ZeLF rect_2x1} % is used to refer this table in the text
\end{table}

The left panel of Fig.~\ref{Fig: fixedpt_2x1} shows these particle equilibria along with simulated particle trajectories. At the origin is an unstable node with the two eigenvalues differing from each other by an order of magnitude (see Table~\ref{table:ZeLF rect_2x1}). This results in the trajectories moving away from this unstable node at different rates in the $r$ and $z$ directions with a larger rate in the $z$ direction. The particle equilibria near the center of the edges of the rectangle are stable nodes. For the stable nodes near the center of the top and bottom edges, the eigenvalue corresponding to the $z$ direction is approximately $30$ times larger in magnitude as compared to the eigenvalue in the $r$ direction. Similarly, for the stable nodes near the left and right edges, the eigenvalue corresponding to the $r$ direction is approximately an order of magnitude larger compared to the eigenvalue in the $z$ direction. Moreover, the saddle points are located near the corners of the rectangle and the corresponding eigenvalues also differ by a factor of $2$. This large disparity in the eigenvalues of the stable nodes and saddle points result in the formation of a slow manifold that connects all these fixed points; specifically consisting of heteroclinic orbits connecting each saddle to the two nearest stable nodes. On this manifold, the eigenvectors corresponding to the smaller magnitude eigenvalues are tangential to the slow manifold. Thus, as can be seen from particle trajectories in the left panel of Fig.~\ref{Fig: fixedpt_2x1}, a typical particle quickly ``snaps" onto the slow manifold and then migrates slowly along this slow manifold towards a stable equilibrium point.

\subsection{Small $\tilde{a}$ limit}

In this section we take the limit $\tilde{a}\xrightarrow{}0$, for which the inertial lift force is negligible and the system of nonlinear ODEs in Eqs.~\eqref{Zelf_r} and \eqref{Zelf_z} reduces to
\begin{equation}\label{large_kappa_2x1}
    \frac{d r_p}{d t}={\frac{1}{2\tilde{R}}D^{2\times1}_r(r_p,z_p)}\:\:\:\text{and}\:\:\:\frac{d z_p}{d t}=\frac{1}{2\tilde{R}}{D^{2\times1}_z(r_p,z_p)}.
\end{equation}
Solving for the fixed points of Eq.~\eqref{large_kappa_2x1} gives us
\begin{align*}
(r^*,z^*)&=\left\{(2,z_p), (-2,z_p) , (r_p,1), (r_p,-1), (0,1/\sqrt{5}), (0,-1/\sqrt{5})\right\}.
\end{align*}
Here the fixed points $(2,z_p), (-2,z_p) , (r_p,1)$ and $(r_p,-1)$ correspond to the boundaries of the square cross-section and hence we ignore them for our linear stability analysis. To determine the stability of the fixed points $(0,\pm 1/\sqrt{5})$, we calculate the eigenvalues of the linear stability matrix in Eq.~\eqref{lin_sys_2x1} giving us
\begin{equation*}
    \lambda_{1,2} = \pm \frac{4i C}{\tilde{R}} \left(\frac{2}{5}\right)^{3/2}.
\end{equation*}
Since the real part of the eigenvalue is zero and the imaginary part is non-zero, the nonlinear nature of the fixed point is inconclusive and it could either be a center or a stable/unstable spiral~\cite{strogatz2019nonlinear}. However, in this simplified secondary drag model, we can analytically solve for the particle trajectories by dividing the two equations in Eq.\eqref{large_kappa_2x1} giving us
\begin{equation*}
\frac{\text{d}r_p}{\text{d}z_p}=\frac{D^{2\times1}_r(r_p,z_p)}{D^{2\times1}_z(r_p,z_p)}=\frac{(1-0.25 r_p^2)(1-5 z_p^2)}{r_p z_p (1-z_p^2)}.    
\end{equation*}
Separating the variables and integrating gives us the following equation for particle trajectories within the cross-section:
\begin{equation}
    z_p (1-z_p^2)^2 (r_p^2-4)^2 = k,
\end{equation}
where $k$ is an arbitrary constant giving different closed curves for the particle trajectories as shown in the right panel of Fig.~\ref{Fig: fixedpt_2x1}. Hence, we conclude the fixed point is a true center. Here, the particle follows a streamline of the secondary flow induced by the cross-sectional vortices resulting in closed trajectories. Moreover, near the fixed point we can get a more simplified equation for closed curves by solving the matrix equation \eqref{lin_sys_2x1}. This gives,
\begin{equation*}
\dot{r}_1=\mp \frac{4 C}{\sqrt{5}\tilde{R}} z_1\,\,\text{and}\,\,\dot{z}_1=\pm \frac{8 C}{25 \sqrt{5}\tilde{R}} r_1,
\end{equation*}
which, upon eliminating the time variable and integrating, results in
\begin{equation*}
r^2_1+\frac{25}{2}z^2_1=c,
\end{equation*}
where $c$ is an arbitrary constant. Thus, the closed curves are concentric ellipses with eccentricity $\sqrt{23}/5$ near the center fixed point.

\begin{figure}
\centering
\includegraphics[width=0.95\columnwidth]{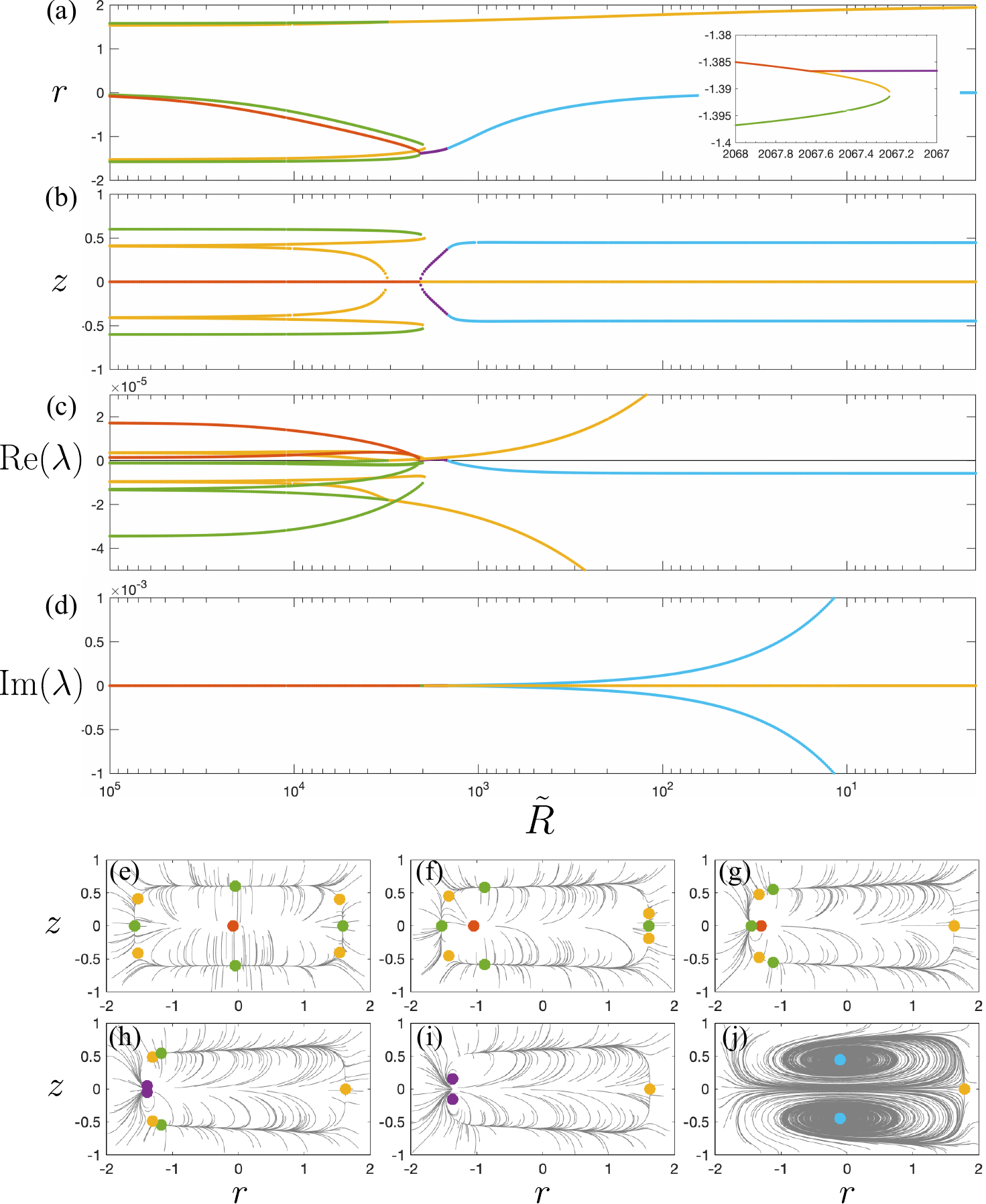}
\caption{Bifurcations in particle equilibria inside a $2\times1$ rectangular cross-section as a function of the dimensionless bend radius $\tilde{R}$ for a fixed dimensionless particle size $\tilde{a}=0.05$. The (a) radial $r$ and (b) vertical $z$ location of the particle equilibria as well as the (c) real and (d) imaginary parts of the eigenvalues $\lambda$ are shown as a function of $\tilde{R}$ (Note that $\tilde{R}$ decreases from left to right). (e)-(j) show the particle equilibria and particle trajectories (gray curves) in the cross-section for $\tilde{R}=10^5$, $3500$, $2200$, $2050$, $1900$ and $100$, respectively. The filled circles denote the fixed points with the size matched with the size of the particle and the color denoting the type of equilibrium point: unstable node (red), stable node (green), saddle point (yellow), unstable spiral (purple) and a stable spiral (cyan). The eigenvalue curves in panels (c) and (d) corresponding to particle equilibria have also been color-coded using the same convention.}
\label{Fig: zelf_2x1}
\end{figure}

\subsection{Bifurcations as a function of $\tilde{R}$}\label{sec: fullzelf}

To investigate the particle equilibria and their bifurcations in the regime where both the inertial lift force and the secondary drag force have non-negligible contribution, we numerically solve for the particle equilibria and the corresponding eigenvalues of the linear stability matrix for the full nonlinear system presented in Eqs.~\eqref{Zelf_r} and \eqref{Zelf_z}.

Figure.~\ref{Fig: zelf_2x1} shows the evolution of particle equilibria and the corresponding eigenvalues as a function of $\tilde{R}$ for a fixed $\tilde{a}=0.05$. At relatively large values of $\tilde{R}$, we first observe a subcritical pitchfork bifurcation where the stable node near the center of the right edge merges with the two saddle points near the top-right and bottom-right corners, leaving behind a single saddle point (see Fig.~\ref{Fig: zelf_2x1}(e)-(g)). As $\tilde{R}$ is further decreased, the stable nodes at the center of the top and bottom edges drifts towards the saddle points to their left and they merge in saddle-node bifurcations (see Fig.~\ref{Fig: zelf_2x1}(g)-(h)). In the same range of $\tilde{R}$ values, the unstable node located near the center of the rectangle first drifts left towards the stable node near the center of the left edge. At first it appears that the unstable node and the stable node merge and give out two unstable spirals on either side in the vertical direction (see Fig.~\ref{Fig: zelf_2x1}(g)-(i)). However, resolving this bifurcation in detail (see the inset in Fig.~\ref{Fig: zelf_2x1}(a)) reveals that the unstable node first undergoes a supercritical pitchfork bifurcation and changes into a saddle point along with two unstable nodes on either side of the saddle point in the vertical direction. The unstable nodes then transition to unstable spirals while the newly created saddle point merges with the stable node on its left in a saddle-node bifurcation. 
%However, resolving this bifurcation numerically, it appears that the stable node first undergoes a subcritical pitchfork bifurcation producing two stable nodes in the vertical direction with a saddle point in between. The saddle point then merges with the unstable node while the two stable nodes change into stable spirals. This is depicted in the inset of Figure~\ref{Fig: rect_2x1}(a) where the saddle-node bifurcation is clearly visible, while the accurate resolution of the transition from stable nodes to stable spirals, which happens over a very narrow region of bend radii, is beyond the scope of the present work. 
As $\tilde{R}$ is further decreased, the two unstable spirals undergo a supercritical Hopf bifurcation and turn into stable spirals (see Fig.~\ref{Fig: zelf_2x1}(i)-(j)). At the smallest bend radius that is physically possible in the theoretical setup, $\tilde{R}_{min}=W/H=2$, two stable spirals are located near the center of the duct while a saddle point is located near the center right edge close to the right wall. This sequence of bifurcations matches well with that for the same particle size in a $2\times1$ rectangular cross-section obtained using the more complete model of Harding et al.~\cite{harding_stokes_bertozzi_2019} as shown by Valani et al.~\cite{Valani2021}. Moreover, Valani et al.~\cite{Valani2021} also investigated large particle sizes and showed that these result in qualitatively different kinds of bifurcations. Hence, the ZeLF model presented here is most accurate for $0<\tilde{a} \lesssim 0.05$.

%We note that from Fig.~\ref{Fig: zelf_2x1}(c) that in the limit $\kappa\xrightarrow{}\infty$, the real part of the eigenvalue for the particle equilibria near the center of the duct seems to approach a non-zero negative constant value suggesting that these particle equilibria remain a stable spiral for large $\kappa$. This is in contrast to what we observed in the large $\kappa$ approximation (inertial lift neglected) where in the right panel of Fig.~\ref{Fig: fixedpt_2x1}, a true center fixed point was realized. This suggests that the non-negligible inertial lift component contributes in the limit of large $\kappa$ corresponding to a vanishingly small particle size and turns the center into a stable spiral. The location of the center and the stable spiral coincide in the limit $\kappa\xrightarrow{}\infty$ but due to the differences in the Jacobian matrices, the nature of the particle equilibria changes. Since the fixed point approaches the location $(r^*,z^*)=(0,\pm1/\sqrt{5})$ as $\kappa\xrightarrow{}\infty$ we have ${\partial D_r}/{\partial r_p}\big|_{(r^*,z^*)}={\partial D_z}/{\partial z_p}\big|_{(r^*,z^*)}\xrightarrow{}0$ while ${\partial L_r}/{\partial r_p}\big|_{(r^*,z^*)}\xrightarrow{}-0.0149$ and ${\partial L_z}/{\partial z_p}\big|_{(r^*,z^*)}\xrightarrow{}-0.7329$ resulting in a negative real part of the eigenvalues and hence a stable spiral.

\section{Particle equilibria and bifurcations in a $1\times2$ rectangular cross-section}\label{sec: 1x2_bifur}

\subsection{Large $\tilde{R}$ limit}

\begin{figure}[t]
\centering
\includegraphics[width=0.5\columnwidth]{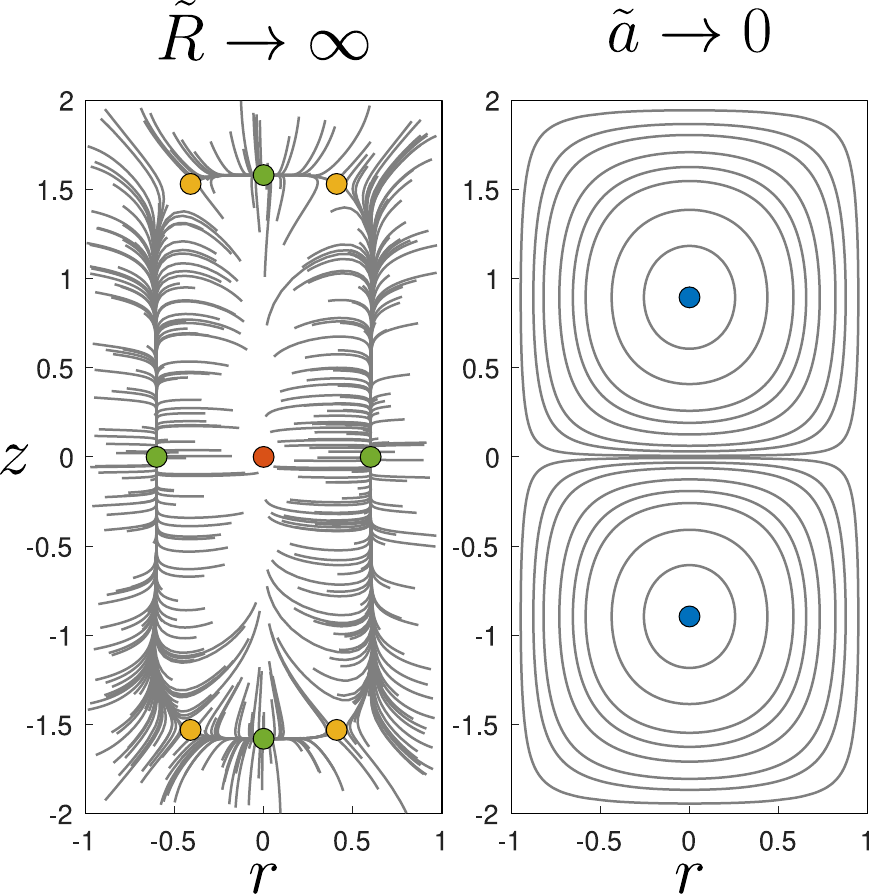}
\caption{Particle equilibria (filled circles) and particle trajectories (gray curves) in the limit of $\tilde{R} \xrightarrow{} \infty$ (left) and $\tilde{a} \xrightarrow{}0$ (right) for a $1\times2$ rectangular cross-section. The color of the filled circles indicates the type of particle equilibria obtained from linear stability analysis: unstable node in red, saddle point in yellow, stable node in green and a center in blue.}
\label{Fig: fixedpt_1x2}
\end{figure}

Since the inertial lift force field in a $1\times2$ rectangular cross-section is modeled as a $90^{\circ}$ rotation of the force field in a $2\times1$ rectangular, the corresponding fixed points also have their $r$ and $z$ co-ordinate switched and are given as follows:
\begin{align*}
(r^*,z^*)=&\big{\{}(0,0), (\pm 0.6000,0) , (0,\pm 1.5800),\\ \nonumber
&(0.4094,\pm 1.5303), (-0.4094,\pm 1.5303)\big{\}}.
\end{align*}
Applying a small perturbation and solving the resulting linear system, we calculate the eigenvalues for each of the nine fixed points. These eigenvalues and the corresponding eigenvectors along with the nature of the fixed points are presented in Table~\ref{table:ZeLF rect_1x2}. Typical particle trajectories along with the fixed points are depicted in the left panel of Fig.~\ref{Fig: fixedpt_1x2}.

\begin{table}
\caption{Eigenvalues $(\lambda_1,\lambda_2)$ and eigenvectors $(\mathbf{V}_1,\mathbf{V}_2)$ for the fixed points $(r^*,z^*)$ of the ZeLF model with a rectangular $1\times2$ cross-section in the limit $\tilde{R}\xrightarrow{}\infty.$} % title of Table
\centering % used for centering table
\begin{tabular}{c c c c c c} % centered columns (4 columns)
\hline\hline %inserts double horizontal lines
Location $(r^*,z^*)$ & Type & $\lambda_1$ & $\lambda_2$ & $\mathbf{V}_1$ & $\mathbf{V}_2$ \\ [0.5ex] % inserts table
%heading
\hline % inserts single horizontal line
\\
(0,0) & unstable node & 0.0110\,$\tilde{a}^3$ & 0.1373\,$\tilde{a}^3$ & $\begin{bmatrix} 0\\ 1 \end{bmatrix}$ & $\begin{bmatrix} 1\\ 0 \end{bmatrix}$ \\[3ex]  % inserting body of the table

($\pm$ 0.6000,0) & stable node & -0.0092\,$\tilde{a}^3$ & -0.2762\,$\tilde{a}^3$ & $\begin{bmatrix} 0\\ 1 \end{bmatrix}$ & $\begin{bmatrix} 1\\ 0 \end{bmatrix}$ \\[3ex] % inserting body of the table 

(0,$\pm$ 1.5800) & stable node & -0.1059\,$\tilde{a}^3$ & -0.0094\,$\tilde{a}^3$ & $\begin{bmatrix} 0\\ 1 \end{bmatrix}$ & $\begin{bmatrix} 1\\ 0 \end{bmatrix}$ \\[3ex] % inserting body of the table 

($\pm$ 0.4092, 1.5304) & saddle & -0.0770\,$\tilde{a}^3$ & 0.0283\,$\tilde{a}^3$ & $\begin{bmatrix} \pm 0.6748\\ 0.7380 \end{bmatrix}$ & $\begin{bmatrix} 0.8526\\ \mp 0.5225 \end{bmatrix}$ \\[3ex] % inserting body of the table

($\pm$ 0.4092, -1.5304) & saddle & -0.0770\,$\tilde{a}^3$ & 0.0283\,$\tilde{a}^3$ & $\begin{bmatrix} \mp 0.6748\\ 0.7380 \end{bmatrix}$ & $\begin{bmatrix} 0.8526\\ \pm 0.5225  \end{bmatrix}$ \\[3ex] % inserting body of the table

%[1ex] % [1ex] adds vertical space
\hline %inserts single line
\end{tabular}
\label{table:ZeLF rect_1x2} % is used to refer this table in the text
\end{table}

\subsection{Vanishingly small particle}

Taking the limit $\tilde{a}\xrightarrow{}0$ in Eqs.~\eqref{Zelf_r} and \eqref{Zelf_z} and using the secondary drag force field for a rectangular $1\times2$ cross-section presented in Sec.~\ref{Sec: zelf 1x2}, we get the following fixed points
\begin{align*}
(r^*,z^*)&=\left\{(1,z_p), (-1,z_p) , (r_p,2), (r_p,-2), (0,2/\sqrt{5}), (0,-2/\sqrt{5})\right\}.
\end{align*}
Here the fixed points $(1,z_p), (-1,z_p) , (r_p,2)$ and $(r_p,-2)$ correspond to the boundaries of the rectangular cross-section and hence we ignore them for our linear stability analysis. Applying a small perturbation and performing a linear stability analysis gives us the following eigenvalues for the fixed points $(0,\pm 2/\sqrt{5})$,
\begin{equation*}
    \lambda_{1,2} = \pm \frac{8i C}{\tilde{R}} \left(\frac{2}{5}\right)^{3/2}.
\end{equation*}
Performing a similar analysis as it was done for a $2\times1$ cross-section, we analytically solve for the particle trajectories by dividing Eqs.~\eqref{Zelf_r} and \eqref{Zelf_z} in the limit $\tilde{a}\xrightarrow{}0$, giving us:
\begin{equation*}
\frac{\text{d}r_p}{\text{d}z_p}=\frac{D^{1\times2}_r(r_p,z_p)}{D^{1\times2}_z(r_p,z_p)}=\frac{(1-r_p^2)(1-1.25 z_p^2)}{4 r_p z_p (1-0.25 z_p^2)}.    
\end{equation*}
Separating the variables and integrating gives us the following equation for particle trajectories within the cross-section:
\begin{equation*}
    z_p (z_p-2)^2 (z_p+2)^2 (1-r_p^2)^2 = k,
\end{equation*}
where $k$ is an arbitary constant giving different closed curves around the fixed point (see right panel of Fig.~\ref{Fig: fixedpt_1x2}). Moreover, near the fixed point we can derive a more simplified equation for closed curves giving us,
\begin{equation*}
\dot{r}_1=\mp \frac{2 C}{\sqrt{5}\tilde{R}} z_1\,\,\text{and}\,\,\dot{z}_1=\pm \frac{64 C}{25 \sqrt{5}\tilde{R}} r_1,
\end{equation*}
which, on eliminating time, gives
\begin{equation*}
r^2_1+\frac{25}{32}z^2_1=c
\end{equation*}
where $c$ is an arbitrary constant. Hence we see that these are equations of concentric ellipses with eccentricity $\sqrt{14}/8$ that correspond to the streamlines of the secondary flow.

\subsection{Bifurcations as a function of $\tilde{R}$}

\begin{figure}
\centering
\includegraphics[width=\columnwidth]{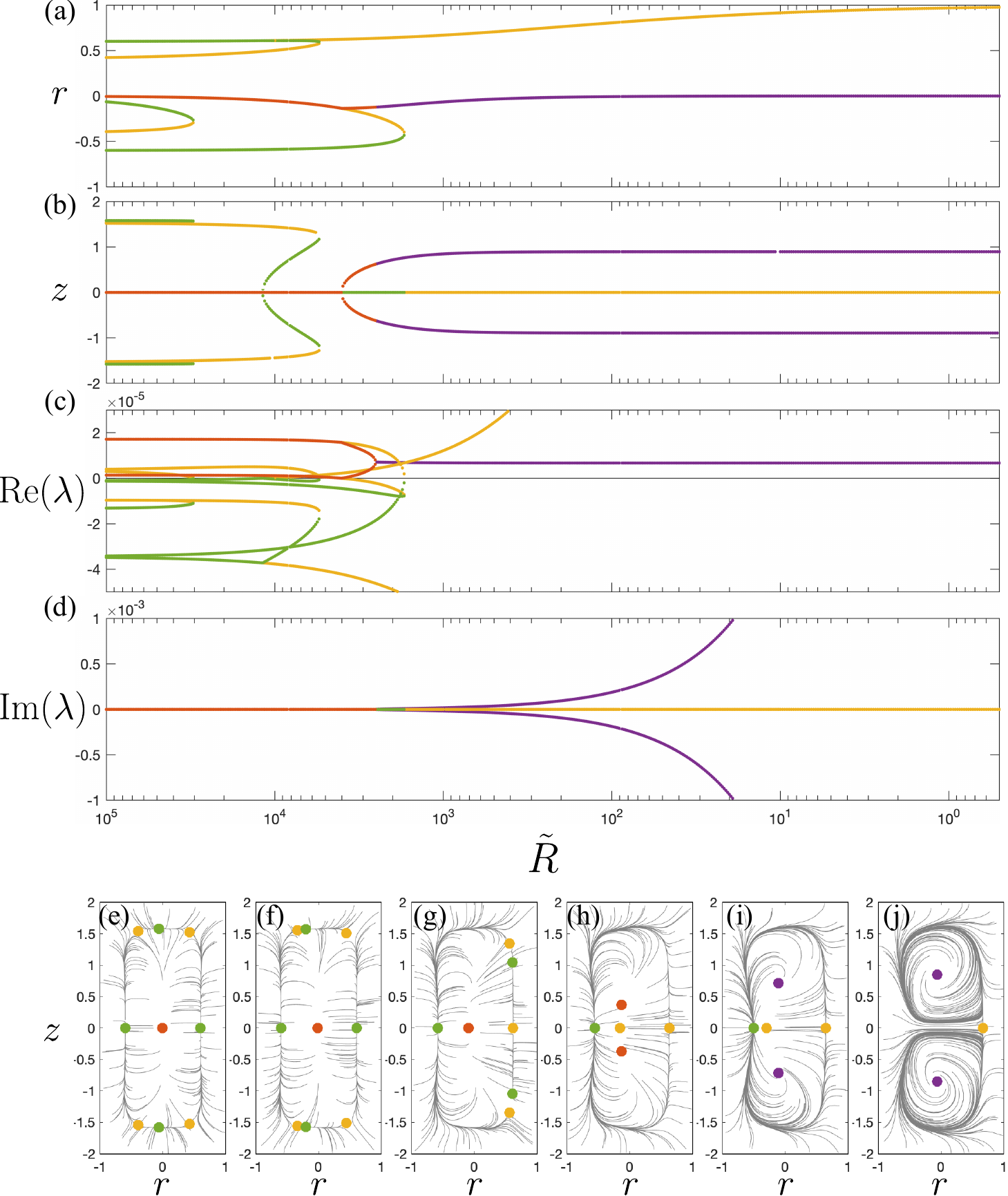}
\caption{Bifurcations in particle equilibria inside a $1\times2$ rectangular cross-section as a function of the dimensionless bend radius $\tilde{R}$ for a fixed dimensionless particle size $\tilde{a}=0.05$. The (a) radial $r$ and (b) vertical $z$ location of the particle equilibria as well as the (c) real and (d) imaginary parts of the eigenvalues $\lambda$ are shown as a function of $\tilde{R}$ (Note that $\tilde{R}$ decreases from left to right). (e)-(j) show the particle equilibria and particle trajectories (gray curves) in the cross-section for $\tilde{R}=10^5$, $35000$, $6000$, $3500$, $2000$ and $1000$, respectively. The filled circles denote the fixed points with the size matched with the size of the particle and the color denoting the type of equilibrium point: unstable node (red), stable node (green), saddle point (yellow), unstable spiral (purple) and a stable spiral (cyan). The eigenvalue curves in panels (c) and (d) corresponding to particle equilibria have also been color-coded using the same convention.}
\label{Fig: zelf_1x2}
\end{figure}

The evolution of particle equilibria in a $1\times2$ rectangular cross-section as a function of $\tilde{R}$ for a fixed $\tilde{a}=0.05$ is depicted in Fig.~\ref{Fig: zelf_1x2}. We first observe saddle-node bifurcations between the two stable nodes near the center of top and bottom edges and the saddle points on their left, at relatively large bend radii. As the bend radius is decreased, the stable node at the center of the right edge undergoes a supercritical pitchfork bifurcation and turns into a saddle point releasing two additional stable nodes on either side in the vertical direction (see Fig.~\ref{Fig: zelf_1x2}(e)-(f)). These newly formed stable nodes undergo saddle-node bifurcations with the two saddle points adjacent to them (see Fig.~\ref{Fig: zelf_1x2}(g)-(h)). As the bend radius is further decreased, the unstable node at the center migrates left and undergoes a supercritical pitchfork bifurcation where it turns into a saddle point releasing two unstable nodes vertically, one on either side (see Fig.~\ref{Fig: zelf_1x2}(h)). The saddle point then goes on to merge with the stable node at the center of the left edge in a saddle-node bifurcation while the two new unstable nodes turn into unstable spirals and develop encompassing limit cycles (see Fig.~\ref{Fig: zelf_1x2}(h)-(j)). Thus, at the smallest physically possible bend radii of $\tilde{R}_{min}=W/W=1$, we have three fixed points: a pair of unstable spirals with limit cycles surrounding them and a saddle point on the right side. This sequence of bifurcations also matches well with that for the same particle size in a $1\times2$ rectangular cross-section obtained using the more complete model of Harding et al.~\cite{harding_stokes_bertozzi_2019} as shown by Valani et al.~\cite{Valani2021}. However, similar to the $2\times1$ cross-section, the ZeLF model for the $1\times2$ cross-section presented here is most accurate for $0<\tilde{a} \lesssim 0.05$ with the bifurcations for bigger particle sizes being qualitatively different and not captured by this model~\cite{Valani2021}.

\section{Conclusions}\label{sec: concl}

We have presented the bifurcations in particle equilibrium positions as given by a reduced-order model of the dynamics of a particle suspended in a fluid flow through a curved duct with a $2\times1$ and a $1\times2$ rectangular cross-section. In the limit of $\tilde{R}\xrightarrow{}\infty$ or $\tilde{a}\xrightarrow{}0$, we analytically obtained the fixed points and their eigenvalues. In the very large bend radius limit where secondary drag is neglected, we obtained stable nodes, saddle-points and an unstable node. Moreover, the large disparity in the magnitude of the two eigenvalues of saddle points and stable nodes resulted in the emergence of a slow manifold for particle dynamics. For a very small particle size where inertial lift force is neglected, we were able to completely solve for the particle dynamics and obtained closed trajectories for particles within the cross-section. Exploring the regime of finite bend radii for a small but non-zero particle size, a number of different bifurcations were observed including saddle-node, pitchfork and Hopf bifurcations. These bifurcations match well with the simulations from the leading order force model as shown in Valani et al.~\cite{Valani2021}. The present simplified ZeLF model is only valid for relatively small particles ($0<\tilde{a}\lesssim 0.05$) and breaks down for larger particle where qualitatively different types of bifurcations are observed~\cite{Valani2021}. Hence, one future direction of the present work would be to formulate a generalized ZeLF model that can accurately capture the particle dynamics and the bifurcations in the particle equilibria for a wider range of particle sizes.
\newline\newline
\textbf{Acknowledgements.} This research is supported under Australian Research Council’s Discovery Projects funding scheme (project number DP160102021 and DP200100834). The results were computed using supercomputing resources provided by the Phoenix HPC service at the University of Adelaide and the Raapoi HPC service at Victoria University of Wellington.

%\bibliography{references}
%\bibliographystyle{unsrtnat}
\printbibliography %Prints bibliography

%\begin{thebibliography}{6}
%

%\bibitem {smit:wat}
%Smith, T.F., Waterman, M.S.: Identification of common molecular subsequences.
%J. Mol. Biol. 147, 195?197 (1981). \url{doi:10.1016/0022-2836(81)90087-5}
%
%\bibitem {may:ehr:stein}
%May, P., Ehrlich, H.-C., Steinke, T.: ZIB structure prediction pipeline:
%composing a complex biological workflow through web services.
%In: Nagel, W.E., Walter, W.V., Lehner, W. (eds.) Euro-Par 2006.
%LNCS, vol. 4128, pp. 1148?1158. Springer, Heidelberg (2006).
%\url{doi:10.1007/11823285_121}%
%
%\bibitem {fost:kes}
%Foster, I., Kesselman, C.: The Grid: Blueprint for a New Computing Infrastructure.
%Morgan Kaufmann, San Francisco (1999)
%
%\bibitem {czaj:fitz}
%Czajkowski, K., Fitzgerald, S., Foster, I., Kesselman, C.: Grid information services
%for distributed resource sharing. In: 10th IEEE International Symposium
%on High Performance Distributed Computing, pp. 181?184. IEEE Press, New York (2001).
%\url{doi: 10.1109/HPDC.2001.945188}
%
%\bibitem {fo:kes:nic:tue}
%Foster, I., Kesselman, C., Nick, J., Tuecke, S.: The physiology of the grid: an open grid services architecture for distributed systems integration. Technical report, Global Grid
%Forum (2002)
%
%\bibitem {onlyurl}
%National Center for Biotechnology Information. \url{http://www.ncbi.nlm.nih.gov}
%
%
%\end{thebibliography}
\end{document}